%% file: main.tex
\documentclass[11pt,reqno]{amsart}
\usepackage[T1]{fontenc}
\usepackage{lmodern}
\usepackage{amsmath,amssymb,amsthm,mathtools}
\usepackage[a4paper,margin=29mm]{geometry}
\usepackage{microtype}
\usepackage[hidelinks]{hyperref}
\newtheorem{theorem}{Theorem}[section]
\newtheorem{lemma}[theorem]{Lemma}
\newtheorem{proposition}[theorem]{Proposition}
\newtheorem{corollary}[theorem]{Corollary}
\theoremstyle{definition}

\theoremstyle{remark}
\newtheorem{remark}[theorem]{Remark}
\newcommand{\F}{\mathbb F_3}

\title[Affine extensions and central compression]{Affine extensions and central compression\newline for product-one-free sequences}
\author{Andreas Volkmann}
\date{22 October 2026}
\begin{document}
\begin{abstract}
Let $H_s=H_{27}\times C_3^s$, where $H_{27}$ is the nonabelian group
of order $27$ and exponent three. We study the elements which can be
adjoined to a product-one-free sequence $L$ of length $2s+5$.
Let $Q$ be its central subsequence and let $h$ be the dimension of
the subgroup generated by $Q$. If $2h-|Q|\le1$, all admissible
individual extensions lie in one affine hyperplane of the full
additive group coordinates, avoiding the identity. The proof combines
coefficient identities with finite, explicitly verified configurations
and a transfer valid in every dimension. We also prove an exact
central compression theorem: a central zero-sum-free block of length
$2h-1$ can be replaced by one term after quotienting by an
$(h-1)$-dimensional central subgroup, with a forward implication for
product-one-freeness and admissible extensions.
For the extraspecial group $E_2$ of order $3^5$ and exponent three,
a finite stress obstruction and central compression yield a further
uniform restriction: in any product-one-free sequence of length
$2r+11$ over $E_2\times C_3^r$, the central subsequence is nonempty
and has length at most $2h-2$, where $h$ is its full central rank.
\end{abstract}
\maketitle

\section{Introduction}

For a finite group $G$, a sequence is product-one-free if no nonempty
subsequence has an ordering with product the identity. Its maximal
possible length is the small Davenport constant $\mathsf d(G)$.
Besides determining this constant, one can ask which elements extend
a sequence whose length is close to the maximum.

In an elementary abelian $p$-group, the admissible extensions of a
zero-sum-free sequence of length $\mathsf d(G)-1$ lie in a proper
affine hyperplane. This is part of the classical $p$-group case of
the invariant $\nu(G)$; see \cite[Theorem~5.5.9]{ghk2006} and
\cite[Theorem~3.5]{geroldinger-yang2026}. We use its familiar
augmentation-ideal proof as background, and investigate a
nonabelian counterpart with explicit hypotheses on the central terms.

Put $H_s=H_{27}\times C_3^s$, where
$H_{27}=\mathrm{UT}_3(\mathbb F_3)$.
One has $\mathsf d(H_s)=2s+6$ for every $s\ge0$
\cite{volkmann2026heisenberg}. We use the full additive coordinates
$W\cong\mathbb F_3^{s+3}$ underlying the class-two group law.
These include the commutator coordinate, not only the abelianization.
For a product-one-free sequence $L$, let
\[
 \mathcal E(L)=\{g\in H_s:L\mathbin{\cdot}g
                  \text{ is product-one-free}\}.
\]
The central subsequence $Q$ contains every position of $L$ whose
value is in $Z(H_s)$, and
$h=\dim_{\mathbb F_3}\langle Q\rangle$ is its full central rank.
Since $Q$ is zero-sum-free, $|Q|\le2h$.

\begin{theorem}[Affine extensions]\label{thm:main-extensions}
Let $s\ge0$ and let $L$ be a product-one-free sequence over $H_s$
of length $2s+5$. If $2h-|Q|\le1$, there exists a linear functional
$\nu:W\to\mathbb F_3$ such that
\[
 \nu(g)=1\qquad(g\in\mathcal E(L)).
\]
If $\mathcal E(L)$ is nonempty, $\nu$ is nonzero.
\end{theorem}

The theorem applies in every dimension and includes an empty central
subsequence. It is a containment theorem, not a classification of
the sequence or of every point in its extension set. The functional
need not vanish on the commutator coordinate; accordingly its level
set is an affine hyperplane in $W$, and is not asserted to be a coset
of a subgroup of $H_s$.

The second ingredient is an exact compression of central blocks.
For a zero-sum-free central block $Q_0$ of length $2h-1$ spanning
$Z_0\cong C_3^h$, put $w=\sum Q_0$. We prove that there is a
hyperplane $U\le Z_0$ with
\[
 \prod_{q\in Q_0}(1-X^q)=N_U(1-X^w).
\]
Quotienting the ambient group by $U$ and replacing $Q_0$ by the
single image of $w$ sends product-one-free sequences to
product-one-free sequences. The same implication holds after any
one admissible extension. Only this forward implication is needed;
the converse is false. This transfer reduces the central deficiencies
zero and one in Theorem~\ref{thm:main-extensions} to the cases of
zero or one central term, or to an abelian quotient.

The proof retains indexed occurrences and actual products throughout.
The underlying affine coefficient identity is standard and also
appears in \cite[Lemma~3.1]{volkmann2026heisenberg}. Our use of it
concerns the transfer of actual product-one-free sequences and the
common extension constraint in the nonabelian group.

We also apply central compression to critical sequences over
$G_r=E_2\times C_3^r$, where $E_2$ has order $3^5$ and exponent three.
The cases $\mathsf d(G_r)=2r+10$ for $0\le r\le3$, and a uniform
exclusion of critical sequences without central terms, were proved
in \cite{volkmann2026smallr}. Here a finite stress obstruction gives
the exclusion of exactly one central term for every $r$. It combines
with compression to give the following statement.

\begin{theorem}[Central deficit in a critical sequence]\label{thm:main-deficit}
Let $S$ be a product-one-free sequence of length $2r+11$ over
$E_2\times C_3^r$, where $r\ge0$. Let $Q$ be its entire central
subsequence and put $h=\dim_{\mathbb F_3}\langle Q\rangle$.
Then $Q$ is nonempty and $|Q|\le2h-2$.
\end{theorem}

Both the one-central-term input and its transfer from the finite
symplectic geometry are given in Section~\ref{sec:q1}. The theorem
follows from Theorem~\ref{thm:q01-exclusion} and
Corollary~\ref{cor:central-deficit}.
The affine extension theorem is proved by
Theorem~\ref{thm:fh-small-centre} and
Corollary~\ref{cor:fh-central-deficit}.
These results do not assert the formula $\mathsf d(G_r)=2r+10$ for
arbitrary $r$, or the extension theorem for central deficit at least two.

\input{sections/foundations}
\input{sections/extensions}
\input{sections/q1}
\input{sections/compression}

\section*{Use of AI-assisted tools}
AI-assisted tools were used in developing arguments, checking finite
configurations, searching the literature, and preparing this manuscript.
Such checks do not constitute independent human peer review. All
mathematical assertions used in the proofs are stated explicitly, and
the source code for the finite verification steps is included with
the submission.

\end{document}

%% file: sections/foundations.tex
\section{Coordinates and elementary coefficient identities}\label{sec:foundations}

All vector spaces and coefficients in this paper are over $\mathbb F_3$.
A sequence is an indexed list; equal values in different positions remain
distinct. A sequence is \emph{product-one-free} if no nonempty selection of
its positions can be ordered to have product the identity. We write
$\mathcal E(L)=\{g:L\mathbin{\cdot}g\text{ is product-one-free}\}$.
The small Davenport constant $\mathsf d(G)$ is the largest length of a
product-one-free sequence over $G$.

Write $H_s=H_{27}\times C_3^s$. We use the underlying vector space
\[
 W=E\oplus\mathbb F_3z\oplus R,
 \qquad E=\mathbb F_3^2,\quad R=\mathbb F_3^s,
\]
with multiplication
\[
 x*y=x+y+\tfrac12\omega(x,y)z,
 \qquad \omega(x,y)=\det(x_E,y_E).
\]
Here $1/2=2$ in $\mathbb F_3$. The center is
$Z(H_s)=\mathbb F_3z\oplus R$ and the commutator subgroup is
$\mathbb F_3z$. We also use $V=W/\mathbb F_3z$ and write $v(x)$ for
the quotient coordinate. A sum without an asterisk is always an additive
sum in the indicated vector space.

The commutative group algebra of an additive space $A$ has basis
$X^a$ with $X^aX^b=X^{a+b}$. Its norm element is
$N_A=\sum_{a\in A}X^a$. The noncommutative group algebra of a group is
distinguished explicitly whenever it occurs.

\begin{lemma}[Top coefficient degrees]\label{lem:affine}
Let $A\cong\mathbb F_3^d$, with $d\ge1$, and let $J$ be the augmentation ideal of
$\mathbb F_3[A]$. Then $J^{2d+1}=0$, every element of $J^{2d}$ is a
scalar multiple of $N_A$, and the coefficient function of every element
of $J^{2d-1}$ is affine on $A$.
\end{lemma}
\begin{proof}
Choose a basis $e_1,\ldots,e_d$ and put $T_i=X^{e_i}-1$. Then
\[
 \mathbb F_3[A]=\mathbb F_3[T_1,\ldots,T_d]/(T_1^3,\ldots,T_d^3),
 \qquad J=(T_1,\ldots,T_d).
\]
The largest nonzero total degree is $2d$. Its monomial
$\prod_iT_i^2$ equals $N_A$. In degree at least $2d-1$ the other
monomials are $T_i\prod_{j\ne i}T_j^2$. The coefficient function of
$T_i$ on the three coordinates $0,1,2$ is $(-1,1,0)$, namely
$-1-a_i$, while that of $T_j^2$ is constant one. These monomials
therefore have affine coefficient functions.
\end{proof}

\begin{lemma}[Abelian bound and extensions]\label{lem:abelian-extensions}
For $A\cong C_3^d$, with $d\ge0$, one has $\mathsf d(A)=2d$. If $d\ge1$ and $L$ is a
zero-sum-free sequence of length $2d-1$, its individual zero-sum-free
extensions all lie in an affine hyperplane avoiding zero.
\end{lemma}
\begin{proof}
For $d=0$ the group is trivial and the bound is immediate. Assume $d\ge1$.
The product $F_L=\prod_{x\in L}(1-X^x)$ has coefficient one at zero
when $L$ is zero-sum-free. Thus the nilpotence assertion in
Lemma~\ref{lem:affine} excludes length $2d+1$. Two copies of each
basis vector attain length $2d$.
For length $2d-1$, write the coefficient function of $F_L$ as
$1+\lambda$ with $\lambda$ linear. If $g$ is an admissible extension,
there is no selection of $L$ with sum $-g$ (including the empty
selection). Hence $1+\lambda(-g)=0$ and $\lambda(g)=1$.
If there are no extensions, the containment is vacuous.
\end{proof}

\begin{lemma}[Supports with central product]\label{lem:support-types}
Let $A$ be a nonempty selection of a sequence over $H_s$ with
$\sum_{i\in A}v(x_i)=0$. Form its graph by joining positions whose
values do not commute. Exactly one of the following types occurs.
\begin{enumerate}
\item The graph has no edges. All orders have the single product
$\sum_{i\in A}x_i$.
\item The graph consists of one triangle and isolated vertices.
The three triangle projections in $E$ are $a,b,-a-b$, with
$\det(a,b)\ne0$, and every isolated position is central in $H_s$.
The products over all orders are exactly
$\sum_{i\in A}x_i+z$ and $\sum_{i\in A}x_i-z$.
\item The products over all orders are the three elements
$\sum_{i\in A}x_i+\mathbb F_3z$.
\end{enumerate}
For a product-one-free sequence the third type is impossible.
In the second type the additive sum in $W$ must be zero.
\end{lemma}
\begin{proof}
The product of an ordering $i_1,\ldots,i_k$ is
\[
 \sum_{i\in A}x_i+\tfrac12
 \sum_{p<q}\omega(x_{i_p},x_{i_q})z.
\]
Two disjoint edges can be made into two adjacent pairs. Independently
swapping those pairs changes the correction by two nonzero scalars.
The four resulting corrections cover $\mathbb F_3$, proving the
third alternative whenever the graph has two disjoint edges.

A graph with an edge but no two disjoint edges is a star or a triangle,
together with isolated vertices. No vertex can have degree one here:
if $j$ were the only neighbor of $i$, the equality
$\omega(x_i,\sum_{a\in A}x_a)=0$ would contradict
$\omega(x_i,x_j)\ne0$. Thus the only remaining nonempty graph is a
triangle. Two triangle projections form a basis of $E$. An isolated
projection is orthogonal to both, hence is zero. The three triangle
projections consequently sum to zero. Their three cyclic pairings
have the same nonzero value. Directly evaluating the six orders gives
the two nonzero corrections, which proves the second alternative.
The first alternative is immediate.

Finally $\sum v(x_i)=0$ means that the additive sum is in
$\mathbb F_3z$. The third product set contains the identity. The
second avoids the identity precisely when its additive sum is zero.
\end{proof}

We call the first, second and third alternatives respectively
commuting, triangular and full supports. The triangle of a triangular
support is its \emph{core}. These descriptions concern original
positions and their actual group values.

%% file: sections/extensions.tex
\section{A common hyperplane when there is at most one central term}
\label{sec:extensions}

Throughout this section, $L=(u_i)_{i\in I}$ is product-one-free over
$H_s$, $m=s+2$, and $|I|=2m+1$. Write $e_i\in E$ and $v_i\in V$
for its projections, and write $c=z$ for the commutator basis vector
from Section~\ref{sec:foundations}. The height coordinate of $u_i$ is
denoted by $\zeta_i$. Let $Q$ be the indexed central subsequence of $L$.

\begin{theorem}\label{thm:fh-small-centre}
If $|Q|\le1$, there is a linear functional $\nu:W\to\F$ such that
$\nu(g)=1$ for every $g\in\mathcal E(L)$.
\end{theorem}

The proof occupies this section. Its finite ingredients are explicitly
stated below and verified by the accompanying program
\texttt{anc/fh\_verify.py}. They involve only groups of orders $27$,
$81$ and $243$. All reductions to these ingredients are uniform in $s$.

\subsection{Cofactors, moments, and the extension identity}

For an indexed sublist $A$, let $F_A^V=\prod_{i\in A}(1-X^{v_i})$;
define $F_A^W$ similarly using the full $u_i$. These are products in
\emph{additive}, commutative group algebras. By Lemma~\ref{lem:affine},
there are scalars $C_i$, linear forms $A_{ij}:V\to\F$, and a linear
form $\Lambda:W\to\F$ such that
\begin{align}
 F_{L\setminus i}^V&=C_iN_V,\\
 [X^v]F_{L\setminus\{i,j\}}^V&=c_{ij}+A_{ij}(v),\\
 [X^w]F_L^W&=c_0+\Lambda(w).
 \label{eq:fh-cofactors}
\end{align}
Only noncommuting pairs $i,j$ will be needed in the second identity.
For $a\in E^*=E\setminus\{0\}$ put
$C_a=\sum_{i:e_i=a}C_i$. The subscript on $C_a$ thus denotes an
aggregate, whereas $C_i$ always refers to one original position.

\begin{lemma}\label{lem:fh-moments}
With these definitions,
\[
 \sum_i C_iv_i=0,\qquad
 \sum_i C_iu_i=\Lambda(c)c.
\]
\end{lemma}
\begin{proof}
For a linear coordinate $\ell$, the derivation
$D_\ell X^u=\ell(u)X^u$ gives
\[
 D_\ell F_L=-\sum_i\ell(u_i)X^{u_i}F_{L\setminus i}.
\]
In $\F[V]$ the product $F_L^V$ is zero, so the first identity follows
by applying every linear coordinate on $V$.
For the height coordinate apply the derivation in $\F[W]$ and then
project to $\F[V]$. Every coefficient on the left is
\[
 \sum_{\zeta\in\F}\zeta\bigl(c_0+\Lambda(v)+\zeta\Lambda(c)\bigr)
 =-\Lambda(c).
\]
The right side is $-\sum_i\zeta_iC_i$ at every coefficient. This proves
the height identity and hence the assertion in $W$.
\end{proof}

A balanced triangle is a set $\{a,b,-a-b\}\subset E^*$ with
$\det(a,b)\ne0$. There are eight such triangles, namely the affine
lines avoiding zero. Every unordered noncommuting pair belongs to
exactly one of them. Choose weights $w(a,b)\in\F$ on these pairs
whose sum on each triangle is one, and set
\begin{align*}
 B_w(v)&=\sum_{i<j,\,\omega(u_i,u_j)\ne0}w(e_i,e_j)A_{ij}(v),\\
 D_w(e)&=\sum_{a:\det(a,e)\ne0}C_aw(a,e),\\
 D_w(0)&=0,\qquad \nu_w=\Lambda-B_w\circ v.
\end{align*}

\begin{lemma}[Weighted extension identity]\label{lem:fh-weighted}
For every admissible individual extension $g$,
\begin{equation}
 \nu_w(g)=1+D_w(e(g)).\label{eq:fh-weighted}
\end{equation}
Moreover, if $e=e(g)\ne0$, then
\begin{equation}
 C_a=C_{-a-e}\qquad\bigl(\det(a,e)\ne0\bigr).
 \label{eq:fh-reflection}
\end{equation}
\end{lemma}
\begin{proof}
Give a triangular support its sign $(-1)^{|A|}$ and let $b(T)$ be the
sum of these signs in a product-one-free list $T$. For $T=L\cdot g$,
the signed charge through an old noncommuting edge $ij$ is
$A_{ij}(v(g))$: multiply its affine cofactor in
\eqref{eq:fh-cofactors} by $1-X^{v(g)}$ and evaluate at $-v_i-v_j$.
Through a new noncommuting edge $ig$ the charge is $C_i$.
Lemma~\ref{lem:support-types} ensures that these coefficients count
exactly triangular supports. Every such support has one core, whose
three edge weights sum to one. Consequently
$b(T)=B_w(v(g))+D_w(e(g))$.

On the other hand $F_T^W$ has constant coefficient $\Lambda(g)$.
At additive zero only the empty support and triangular supports can
contribute: a nonempty commuting support there would have product
one, full supports are forbidden, and every triangular support has
additive sum zero. Thus $\Lambda(g)=1+b(T)$, proving
\eqref{eq:fh-weighted}. Summing the new-edge charges over positions
with image $a$ gives the same collection of triangular supports as
summing over their third corners, whose image is $-a-e$. This proves
\eqref{eq:fh-reflection}.
\end{proof}

\begin{lemma}[Two exceptional direction configurations]\label{lem:fh-exceptions}
The desired common functional exists unless one of the following occurs:
\begin{enumerate}
\item[$\mathrm{K8}$:] $C_a=t\ne0$ on all eight points of $E^*$, and
every such direction contains an admissible extension;
\item[$\mathrm{K3}$:] for a balanced triangle $T$, $C_a=t\ne0$ on
$T$ and $C_a=0$ elsewhere, and all three directions of $T$ contain
admissible extensions.
\end{enumerate}
\end{lemma}
\begin{proof}
Let $A\subset E^*$ be the directions of admissible extensions. If $A$
lies on a line through zero, give weight one in each balanced triangle
to an edge avoiding that line, and zero to the other two edges. Then
$D_w=0$ on $A$, so \eqref{eq:fh-weighted} suffices, including at central
extensions.

Otherwise choose independent $e,f\in A$. In coordinates $e=(1,0)$,
$f=(0,1)$, the identifications in \eqref{eq:fh-reflection} have precisely
the two connected classes
\[
 T=\{e,f,-e-f\},\qquad
 E^*\setminus T=\{-e,-f,e+f,e-f,-e+f\}.
\]
Thus $C$ is constant, or takes distinct values $A_0$ on $T$ and $B_0$
off $T$. In the latter case its reflection filter is exactly $T$:
reflection about a point of $T$ preserves the two classes, whereas
reflection about any point outside $T$ joins them. This also follows
by checking the representatives $-e$ and $e-f$; permutations of $T$
are transitive on $-T$ and on the remaining two points. Hence $A\subset T$.

If $C_a=t$ is constant and $t=0$, any weights work. For $t\ne0$, the
possible eight values $D_w(a)$ are exactly the vectors with total sum
$2\cdot8t=t$. Indeed a variation of triangle weights with total zero
produces every zero-sum variation on its three vertices: the change
at a vertex is $-t$ times the change on the opposite edge. These local
zero-sum spaces span the full zero-sum space on the eight points,
since the triangle incidence graph is connected. If $A\ne E^*$,
put the required total on a missing point and set $D_w=0$ on $A$.
The remaining constant case is K8.

In the two-colour case with $B_0\ne0$, an auxiliary triangle through
each $e\in T$ meets $T$ only at $e$. Changing an incident weight and
the opposite weight by opposite amounts adjusts $D_w(e)$ arbitrarily,
without changing its values at the other two points of $T$. Thus
$D_w=0$ on $T$ is possible. If $B_0=0$, the possible three values on
$T$ have exactly total $2A_0$. Unless $A=T$, put this total on a missing
point. The only remaining case is K3.
\end{proof}

\subsection{Central matrices and the exclusion of K8}

Let $Z_Q=\operatorname{span}_W(Q)$, and place a bar over values in
$W/Z_Q$. This is only a linear quotient for now. For an extension
$h_e$ and a balanced triple $a,b,e$, define the matrix, indexed by
original positions in classes $a,b$,
\begin{equation}
 K_{ij}(h_e)=-[X^{-u_i-u_j-h_e}]F_Q^W.
 \label{eq:fh-central-matrix}
\end{equation}
Its row sums are $C_i$ and its column sums are $C_j$. Indeed every
triangular support through $i,h_e$ consists of the third corner $j$
and a sublist of the original central list $Q$. Its additive sum is
zero and its sign is minus the sign of that central sublist.
In particular a nonzero entry implies
$\bar u_i+\bar u_j+\bar h_e=0$. The weighted fibre functions
\[
 \mu_a(x)=\sum_{i:e_i=a,\,\bar u_i=x}C_i
\]
therefore satisfy
\begin{equation}
 \mu_b(-x-\bar h_e)=\mu_a(x).\label{eq:fh-mu-reflection}
\end{equation}
A zero fibre weight does not imply absence of original positions.

\begin{lemma}\label{lem:fh-k8}
K8 is impossible if $\dim\operatorname{span}_V(v(Q))\le4$.
\end{lemma}
\begin{proof}
Put $R_Q=\operatorname{span}_V(v(Q))$. In K8 define
$b_a=t^{-1}\sum_{i:e_i=a}C_i\bar u_i$. Summing
\eqref{eq:fh-central-matrix} gives
$b_a+b_b=-\bar h_e$ when $a+b+e=0$ and $a,b$ are independent.

We use the following elementary eight-point observation. If the sum
$\phi(a)+\phi(b)$ for independent $a,b\in E$ depends only on $a+b$,
then $\phi$ is the restriction of an affine map on $E$. To check this,
put $p=\phi(a)$, $q=\phi(b)$ and $r=\phi(a+b)$ for independent $a,b$.
The pair equations give, respectively at
$-a,-b,-a-b,a-b,-a+b$,
\[
 p-q+r,\quad -p+q+r,\quad -p-q,\quad -p-r,\quad -q-r.
\]
These are the affine values with constant term $p+q-r$.

Thus $b_a=M(a)+d$ and $\bar h_e=M(e)+d$. Project to $V/R_Q$.
Lemma~\ref{lem:fh-moments} gives
$0=\sum_{a\in E^*}t b_a=8td=-td$, so the projected $d$ is zero.
The projections of eight chosen extensions $h_e$ therefore span a
space of dimension at most $2+\dim R_Q$. But the differences
$\nu_w-\nu_{w'}$ factor through $V$, and their evaluations on those
eight extensions, by the constant-profile calculation in
Lemma~\ref{lem:fh-exceptions}, span the seven-dimensional zero-sum
space in $\F^8$. Hence $7\le2+\dim R_Q$, a contradiction.
\end{proof}

\subsection{No central terms}

\begin{lemma}\label{lem:fh-qzero}
Theorem~\ref{thm:fh-small-centre} holds when $Q$ is empty.
\end{lemma}
\begin{proof}
By Lemmas~\ref{lem:fh-exceptions} and \ref{lem:fh-k8}, only K3 needs
to be excluded. Choose one admissible extension $h_a,h_b,h_d$ in
each direction of its balanced triangle $T=\{a,b,d\}$, and put
$K=h_a+h_b+h_d$. The weighted identity gives
\begin{equation}
 \nu_w(K)=2t\ne0,\qquad e(K)=0.\label{eq:fh-Knonzero}
\end{equation}

For a full original value $u$, write $n(u)$ for its multiplicity and
$C(u)$ for its cofactor. We have $n(u)\le2$. When $u$ does not commute
with $h_e$, \eqref{eq:fh-central-matrix} with $Q$ empty gives
\[
 C(u)=-n(-u-h_e).
\]
Thus every value with $C(u)\ne0$ has a genuine original reflected
partner $r_e(u)=-u-h_e$ with nonzero cofactor. The weight
$n(u)C(u)$ is preserved by this reflection. A closed even reflection
path in $E$ acts on full values as a central translation $u\mapsto
u+\delta$. If $\delta\ne0$, its three iterates give three distinct
original values with the same $E$-image, additive sum zero, and hence
commuting product one. Therefore such translations are zero on the
set of values with nonzero cofactor.

Outside $T$ the reflection graph is $K_{2,3}$ on $-T$ and the two
remaining points. In coordinates $a=(1,0)$, $b=(0,1)$, $d=(2,2)$,
the cycle
$(1,2),-a,(2,1),-b,(1,2)$ has colours $b,d,a,d$ and acts as
$u\mapsto u+h_a+h_b-2h_d=u+K$. By \eqref{eq:fh-Knonzero} this is
nonzero. Connectivity shows that no original value outside $T$ has
nonzero cofactor.

Put $S_e=\sum_{i:e_i=e}C_i u_i$ and $M=\sum_iC_i u_i$.
The inner reflections give
$S_b+S_d=-t h_a$ and its two cyclic analogues. Thus $M=tK$.
The $V$-part of Lemma~\ref{lem:fh-moments} implies $v(K)=0$.
Take a value $u_a$ with nonzero cofactor. Two even reflection paths
produce actual original values
$u_b=u_a+h_b-h_a$ and $u_d=u_a+h_d-h_a$.
Their images in $E$ are different, and their additive sum is $K$.
They form a triangular $V$-zero support, so
Lemma~\ref{lem:support-types} gives $K=0$, contradicting
\eqref{eq:fh-Knonzero}.
\end{proof}

\subsection{Finite products and a dimension-independent moment argument}

For a list $B$, let $\Pi_0(B)$ contain the products of all ordered
selections of its indexed positions, including the empty product.
We will repeatedly use a normal subgroup $J\cong H_k\le H_s$ and a
decomposition $L=B_0\mathbin{\cdot}R$. If each $h_e$ is individually
admissible and
$U=\bigcup_e\Pi_0(B_0\mathbin{\cdot}h_e)$, then every nonempty
ordered selection from $R$ whose product belongs to $J$ avoids
$U^{-1}$. Otherwise it could be concatenated with a disjoint selection
of $B_0\mathbin{\cdot}h_e$ to give product one in that particular
extension. In particular, if $U=J$ and $|R|>2(s-k)$, the abelian
bound in $H_s/J\cong C_3^{s-k}$ gives a contradiction.

When $U\ne J$, fix any one order on $R$ and let $f_R$ be the
coefficient function of $\prod_{u\in R}(1-u)$ in that order in the
\emph{noncommutative} group algebra. Then
\begin{equation}
 f_R(1)=1,\qquad
 \operatorname{supp}(f_R|_J)\subset
 A:=\{1\}\cup(J\setminus U)^{-1}.
 \label{eq:fh-support-complement}
\end{equation}

\begin{lemma}[Projected moments]\label{lem:fh-local-moments}
Suppose $J$ contains the commutator subgroup and the quotient $H_s/J$
has dimension $d$. If $P$ is a polynomial of degree at most $b$ in
the abelianization coordinates on $J$, and $b+2d<|R|$, then
\[
 \sum_{g\in J} f_R(g)P(v(g))=0.
\]
\end{lemma}
\begin{proof}
Choose complementary radical coordinates $\rho_1,\ldots,\rho_d$.
The sum in question is the signed sum over all selections of $R$ of
\[
 P\left(\sum v(u)\right)
 \prod_{j=1}^d\left(1-\rho_j\left(\sum v(u)\right)^2\right).
\]
It is the result of $|R|$ translation differences applied to a
polynomial of degree at most $b+2d$. Each difference lowers total
degree, so it vanishes. Only the polynomial is projected to $V$;
the support restriction \eqref{eq:fh-support-complement} concerns
the full group values.
\end{proof}

The finite checks below have a particularly direct exhaustive
interpretation. For a fixed indexed list $(b_1,\ldots,b_n)$ define
$P(\varnothing)=\{1\}$ and
\begin{equation}
 P(M)=\bigcup_{i\in M}\{x*b_i:x\in P(M\setminus\{i\})\}.
 \label{eq:fh-dp}
\end{equation}
Induction on $|M|$ proves that this is exactly the set of all products
of all orders on $M$. Thus $\Pi_0(B)=\bigcup_M P(M)$, and $B$ is
product-one-free precisely when $1\notin P(M)$ for every nonempty
$M$. Formula~\eqref{eq:fh-dp}, with full BCH multiplication, is the
algorithm used by the supplied verifier. Equal values keep separate
indices. The polynomial identities used below are checked at every
point of the indicated finite group. This provides reproducible finite
proofs rather than searches in the parameter $s$.

\subsection{One central term: the nonzero quotient branch}

Assume now $Q=(q)$. If $v(q)=0$, all noncentral cofactors $C_i$ vanish
because their projected products contain $1-X^{v(q)}=0$. The weighted
identity then already gives the conclusion. We therefore assume
$v(q)\ne0$, so $q$ and $c$ are independent. K8 is excluded by
Lemma~\ref{lem:fh-k8}. Suppose K3 occurs, choose its three extensions
$h_a,h_b,h_d$, and set $K=h_a+h_b+h_d$. Bars now mean the linear
quotient $W/\langle q\rangle$.

\begin{lemma}\label{lem:fh-K-in-q}
For every such choice of extensions, $K\in\langle q\rangle$.
\end{lemma}
\begin{proof}
Suppose $\bar K\ne0$. The outer four-cycles in
\eqref{eq:fh-mu-reflection} act by translation by $\bar K$.
On a three-point translation orbit, both its total mass and its
first moment vanish. Hence the outer first moments vanish. The inner
reflection equations give
\[
 \overline{\sum_iC_i u_i}=t\bar K.
\]
Lemma~\ref{lem:fh-moments} implies
\begin{equation}
 \bar K=(\Lambda(c)/t)\bar c\ne0.\label{eq:fh-K-height}
\end{equation}

We detail a disjointness consequence. Inside $T$, the reflection
orbits with nonzero fibre weight consist of one possible fixed
three-fibre orbit $x_e=\bar h_e+\bar K$, and six-fibre orbits.
Indeed, from a fibre $x$ over $a$, the even paths produce
$x$, $x+\bar h_b-\bar h_a$, $x+\bar h_d-\bar h_a$; their sum is
$\bar K$. The other parity gives a second such triple unless the
orbit is the fixed one. Choosing one original representative in each
fibre gives a balanced triple whose product can be put in
$\langle q\rangle$, by \eqref{eq:fh-K-height} and the two possible
corrections $\pm c$. Each six-fibre orbit yields two disjoint triples.
Outside $T$, every translation orbit has three distinct fibres over
each of five directions. Representatives in one such direction
commute and have sum in $\langle q\rangle$, giving five disjoint
blocks. If their numbers are $I_3\in\{0,1\}$, $I_6$ and $O$,
respectively, these blocks together with the unused original $q$
form a zero-sum-free list in $\langle q\rangle$. Therefore
\[
 I_3+2I_6+5O\le1.
\]
Since the inner total weight is $t\ne0$, we have $I_3=1$, $I_6=O=0$.

Choose arbitrary original representatives $u_a,u_b,u_d$ of the three
remaining fibres. Their sum is $\alpha q+k c$ with $k\ne0$.
If $\alpha=0$, the balanced triple has an order with product one;
if $\alpha=2$, the same is true after adjoining the unused original
$q$. Consequently their sum is always
\begin{equation}
 U:=u_a+u_b+u_d=q+k c.\label{eq:fh-U}
\end{equation}
Varying one representative shows that each of these fibres contains
exactly one full original value, of multiplicity $n_e\in\{1,2\}$.

We must also control positions whose fibre weights cancel. For a
nonzero matrix entry \eqref{eq:fh-central-matrix}, there is a unique
$\epsilon\in\{0,1\}$ with
$u_i+u_j+h_e=-\epsilon q$. The original positions $i,j,u_e$ are
distinct, since their $E$-images are the three different triangle
points. Their sum has the form $\alpha q+k c$ with the same nonzero
$k$; by the previous product argument $\alpha=1$. Thus, for each
fixed $e$, the value $\epsilon=\epsilon_e$ is the same for every
nonzero entry of that matrix. This is a statement about actual
supports, so it is unaffected by coefficient cancellation.
Applying a nonzero row to one of the three known full values shows
\begin{equation}
 h'_e:=h_e+\epsilon_e q=u_e-U,
 \qquad h'_a+h'_b+h'_d=U.\label{eq:fh-pure-reflections}
\end{equation}
For every original value $u$ noncommuting with $h_e$, the now fixed
central mask gives the exact row identity
\[
 C(u)=-(-1)^{\epsilon_e}n(-u-h'_e).
\]
Since every multiplicity is at most two, any nonzero matrix entry
has both endpoints of nonzero cofactor. Thus values of nonzero
cofactor have genuine reflection partners under $u\mapsto-u-h'_e$;
the reciprocal relation shows that $n(u)C(u)$ is preserved.
Outer four-cycles translate full values
by $U\ne0$, producing three commuting distinct original values with
sum zero; hence there are no outer values of nonzero cofactor.
An additional inner orbit would consist of six values
$u_e+\delta,u_e-\delta$ ($e\in T$), with $\delta\ne0$ central.
Their sum is $2U$, so together with $q$ they are $V$-zero with six
active positions. Lemma~\ref{lem:support-types} forbids this full
support. All noncentral nonzero cofactors are therefore concentrated
on $u_a,u_b,u_d$. The first moment in $V$ gives $C_q=-t$.

For the remainder write $a=u_a$, $b=u_b$, $d=u_d$. The exact rows give
\[
 h_e=-u_f-u_g-\epsilon_e q,\qquad
 t=-n_f n_g(-1)^{\epsilon_e}.
\]
Let $\kappa=1-\epsilon_a-\epsilon_b-\epsilon_d$. Then
$K=\kappa q+k c$ and Lemma~\ref{lem:fh-moments} gives $\Lambda(c)=tk$.
For each separate extension, new triangular supports have precisely
the two other charged corners and the fixed central mask. Their
total sign is $t$, so $b(L\cdot h_e)=b(L)+t$.
The proof of Lemma~\ref{lem:fh-weighted} consequently gives equal
values $\Lambda(h_e)=1+b(L)+t$, and hence $\Lambda(K)=0$.
It follows that
\[
 0=\kappa\Lambda(q)+k\Lambda(c)=\kappa\Lambda(q)+t,
 \qquad \kappa\ne0.
\]
The possible multiplicity and mask patterns, up to permutation, are
therefore exactly
\begin{center}
\begin{tabular}{c|c|c}
$(n_a,n_b,n_d)$&$t$&$(\epsilon_a,\epsilon_b,\epsilon_d)$\\\hline
$(1,1,1)$&$2$&$(0,0,0)$\\
$(1,1,1)$&$1$&$(1,1,1)$\\
$(2,1,1)$&$2$&$(0,1,1)$\\
$(1,2,2)$&$2$&$(0,1,1)$
\end{tabular}
\end{center}
Indeed the equations for $t$ determine the masks from the
multiplicities; the other mixed-mask choices give $\kappa=0$.
The pattern $(2,2,2)$ is forbidden because those six positions
together with $q$ form a full $V$-zero support.

We now exclude all four rows by the following explicit finite data.
For the $(1,2,2)$ row, the five-position block
$(b,b,d,d,-b-d)$ in $L\cdot h_a$ has product set all of
$J=\langle b,d\rangle\cong H_0$. Its complement has length
$2s+1>2s$, contradicting the abelian bound in $H_s/J$.
This product-set assertion is the first check in the verifier.

For the remaining rows use full coordinates $(x,y,u;\zeta)$ in
$J=\langle a,b,q\rangle\cong H_1$ with
\[
 a=(1,0,0;0),\quad b=(0,1,0;0),\quad q=(0,0,1;0),
 \quad d=(-1,-1,1;k),\qquad k\in\{1,2\}.
\]
These are changes of coordinates, not changes of the original
elements. In the $(2,1,1)$ row put $B_0=(q,a,a,b,d)$.
The separate product sets $\Pi_0(B_0\cdot h_e)$ have sizes
$71,67,67$, and their union is all $81$ elements of $J$, for both
$k$. The complement has length $2s>2(s-1)$, again a contradiction.

In the two $(1,1,1)$ rows put $B_0=(q,a,b,d)$ and
$h_e=-u_f-u_g-\epsilon q$, with common $\epsilon=0$ or $1$.
Each of the three separate product sets has size $50$, and their
union has size $68$ for $\epsilon=0$, or $71$ for $\epsilon=1$.
Define $A$ as in \eqref{eq:fh-support-complement}.
Write $A_v=1-X^v$ in the additive algebra of $V_J$. The following
identities hold at every $g\in J$, for both values of $k$. Define
\[
 \begin{aligned}
 d_0(v)&=[X^{-v}]A_{v(q)}A_{v(a)}A_{v(b)},\\
 d_1(v)&=[X^{-v}]A_{v(q)}A_{v(b)}A_{v(d)}.
 \end{aligned}
\]
Then
\begin{equation}\label{eq:fh-nz-duals}
 \begin{aligned}
 \mathbf1_A(g)\bigl(d_0(v(g))-(u^2-u)\bigr)&=\mathbf1_{\{1\}}(g),\\
 \mathbf1_A(g)\bigl(d_1(v(g))-(1+u^2-x-y+yu)\bigr)&=0.
 \end{aligned}
\end{equation}
The first identity applies to $\epsilon=0$, the second to $1$.
These four pointwise identities and the stated product sets are
verified exhaustively by \eqref{eq:fh-dp} and evaluation at all
$81$ points.

Here $|R|=2s+1$ and $H_s/J\cong C_3^{s-1}$, so quadratic moments
vanish by Lemma~\ref{lem:fh-local-moments}. Convolution of $f_R$ with
$d_0$ is exactly the original cofactor $C_d=t$; convolution with
$d_1$ is $C_a=t$. Summing \eqref{eq:fh-nz-duals} against $f_R$ gives
$C_d=1$ when $t=2$, or $C_a=0$ when $t=1$. Both are contradictions.
This proves $\bar K=0$.
\end{proof}

\subsection{One central term: the zero quotient branch}

\begin{lemma}\label{lem:fh-qone-K3}
K3 is impossible when $Q=(q)$.
\end{lemma}
\begin{proof}
We have already dealt with $v(q)=0$. By Lemma~\ref{lem:fh-K-in-q},
every choice of the three extensions satisfies $K=\kappa q$.
Equation~\eqref{eq:fh-Knonzero}, which follows from the weighted
identity in K3 independently of $Q$, gives $\kappa\ne0$ and
$\nu_w(q)\ne0$. Two extensions in the same direction differ by a
multiple of $q$, by varying that choice in $K$. Their $\nu_w$-values
are equal, so they must coincide. There is exactly one full extension
$h_e$ in each direction of $T$.
Choose coordinates with $q$ radical of height zero and, by a linear
height shear, with $h_a,h_b,h_d$ of height zero. This is possible
because $K=\kappa q$.

First note that three disjoint balanced original triples with sums
$\alpha_iq$ are impossible. Each triple has possible products
$\alpha_iq\pm c$. For any two disjoint triples, choosing opposite
corrections shows that $\alpha_i+\alpha_j$ can be neither $0$ nor
$2$ (in the latter case append the original $q$). Thus every pair
must have sum $1$. For three triples this forces all $\alpha_i=2$;
choosing correction $+c$ for all three then gives product one.

In $W/\langle q\rangle$ the outer reflection graph is the bipartite
$K_{2,3}$. Its cycles now have zero translation, so every orbit of
nonzero fibre weights has five fibres, one per outer direction,
with equal weight. Its three fibres over $-T$ have sum zero and
yield a balanced original triple with sum in $\langle q\rangle$.
Since all outer aggregate weights are zero, any such orbit would
require a second. There is also an inner nonzero orbit because
$t\ne0$. These would yield three disjoint triples, a contradiction.
Hence all outer fibre weights vanish.

The inner orbits are the fixed three-fibre orbit $x_e=\bar h_e$ and
six-fibre orbits. Every six-fibre orbit splits into two disjoint
three-fibre groups with sum zero in the quotient, by the even-path
formula used above. The preceding disjoint-triples argument leaves
exactly two possibilities: the fixed three-fibre orbit alone, or one
six-fibre orbit alone. Positions in zero-weight fibres are retained
in all subsequent complements.

\emph{The six-fibre orbit.}
The two three-fibre groups have sums $\alpha q,\beta q$ with
$\alpha+\beta=1$, for every choice of representatives. Varying one
representative proves that each fibre is monochromatic: it contains
one full original value of multiplicity one or two. Write them as
\[
 u_e^+=h_e+\delta+p_eq,\qquad
 u_e^-=h_e-\delta+m_eq,\qquad \bar\delta\ne0.
\]
In the individually admissible $L\cdot h_e$, the three values
$u_e^+,u_e^-,h_e$ commute and sum to $(p_e+m_e)q$.
Using also the unused original $q$ forces $p_e+m_e=1$.
The total of the six original values is $q$, but also $2K$;
therefore $\kappa=2$.

Each of these monochromatic fibres has nonzero individual cofactor.
A row of \eqref{eq:fh-central-matrix} starting at its plus value can
reach only the minus fibre in the other direction, and conversely.
The full sums with $h_e$ in the two orientations are
$(p_f-p_g)q$ and $(p_g-p_f)q$. A nonzero entry requires each to be
$0$ or $2q$, so $p_f=p_g$. Absorb the common tag in $\delta$ to get
\begin{equation}
 u_e^+=h_e+\delta,\qquad u_e^-=h_e-\delta+q.
 \label{eq:fh-six-values}
\end{equation}
All these cross-entries use the empty central mask. Their row sums
make the three plus multiplicities equal, and the three minus
multiplicities equal. Multiplicity two on either side would give
three disjoint balanced triples. Thus all six multiplicities are one,
their $C_i$ are $2$, and $t=1$.

We also have $v(\delta)\notin\langle v(q)\rangle$. Otherwise write
$v(\delta)=\mu v(q)$. Replacing a plus value in the six-position
support by its admissible $h_e$ produces a support with six active
positions and $V$-sum $(1-\mu)v(q)$. This must equal $v(q)$:
sum zero is full, and sum $2v(q)$ becomes full after adjoining $q$.
Thus $\mu=0$. Replacing a minus value instead gives $\mu v(q)$,
forcing $\mu=1$, a contradiction. These are two separate uses of
individual admissibility.

Consequently $J=\langle h_a,h_b,q,\delta\rangle\cong H_2$ is normal.
In full coordinates $(x,y,u,v;\zeta)$ choose
\[
 \begin{gathered}
 h_a=(1,0,0,0;0),\quad h_b=(0,1,0,0;0),\quad
 h_d=(2,2,2,0;0),\\
 q=(0,0,1,0;0),\qquad \delta=(0,0,0,1;0).
 \end{gathered}
\]
Put $B_0=(q,u_e^+,u_e^-:e\in T)$. Each separate
$B_0\cdot h_e$ is product-one-free, has length eight, and its product
set has $241$ elements. Their union omits exactly
$-\delta$ and $-q+\delta$. These statements are verified at all
indexed masks by \eqref{eq:fh-dp} in the $243$-element group.
Therefore the complement $R$ has length $2s-2$ and its restricted
coefficient function is supported in $\{1,\delta,q-\delta\}$.
The polynomial $P=1-2u-v$ equals one at the identity and zero at
the other two values. Its moment against $f_R$ is thus one.
But $1+2(s-2)<2s-2$, so Lemma~\ref{lem:fh-local-moments} makes
that moment zero. This excludes the six-fibre case.

\emph{The fixed three-fibre orbit.}
Its original values are $h_e+\lambda q$, $\lambda\in\F$.
Two original positions in the same fibre, together with the separate
extension $h_e$, commute; product-one-freeness with the original
$q$ forces their tags to satisfy $\lambda_i+\lambda_j=1$.
Three original positions would therefore all have tag $2$ and
already have product one. The nonempty multiplicity profiles in
each fibre, listed at tags $0,1,2$, are exactly
\begin{equation}
 (1,0,0),\quad(0,1,0),\quad(0,0,1),\quad(1,1,0),\quad(0,0,2).
 \label{eq:fh-five-profiles}
\end{equation}
This counts all original positions in the three fibres, even when
an individual cofactor is zero.

Let $n_e(\lambda)$ denote these counts. Formula
\eqref{eq:fh-central-matrix} gives, at every occupied tag,
\begin{equation}
 C_e(\lambda)=n_f(2-\kappa-\lambda)-n_f(-\kappa-\lambda)
 \label{eq:fh-tag-rows}
\end{equation}
for both other directions $f$. The two values must coincide, and
$\sum_\lambda n_e(\lambda)C_e(\lambda)=t\ne0$ must be the same in
all three directions.

There are $2\cdot5^3=250$ profiles in
\eqref{eq:fh-five-profiles} with $\kappa\in\{1,2\}$.
Equation~\eqref{eq:fh-tag-rows} leaves exactly $26$; requiring each
of the three separate augmented original kernels to be
product-one-free leaves exactly $10$. These are exhaustive finite
checks by \eqref{eq:fh-dp}, not assumptions on other positions.
Six surviving kernels have length six. Up to permuting directions
their profiles are
\[
 \begin{array}{c|c|c}
 \kappa&t&(n_a,n_b,n_d)\\\hline
 1&1&((1,0,0),(0,0,2),(0,0,2))\\
 2&2&((0,1,0),(0,0,2),(0,0,2)).
 \end{array}
\]
For each, every separate augmentation covers all $81$ elements of
$J=\langle h_a,h_b,q\rangle\cong H_1$. Its complement length
$2s-1>2(s-1)$ excludes it by the abelian quotient bound.

The remaining four kernels have length four and are
$B_0=(q,u_a,u_b,u_d)$, where $u_e=h_e+\eta q$, with
\[
 (\kappa,\eta,t)=(1,1,2),(1,2,1),(2,0,1),(2,2,2).
\]
Use coordinates $q=(0,0,1;0)$, $h_a=(1,0,0;0)$,
$h_b=(0,1,0;0)$, $h_d=(2,2,\kappa;0)$ in $J\cong H_1$.
The allowed sets $A$ in \eqref{eq:fh-support-complement} have
sizes $11,11,19,11$, respectively. Put
\[
 \lambda_a(g)=[X^{-v(g)}]
 (1-X^{v(q)})(1-X^{v(u_b)})(1-X^{v(u_d)}).
\]
The following identities hold on the full allowed set $A$:
\begin{center}
\begin{tabular}{c|c|l}
$(\kappa,\eta,t)$&identity on $A$&$P(x,y,u)$\\\hline
$(1,1,2)$&$\mathbf1_{\{1\}}+P=0$&$2+u^2+2yu+y^2+x^2$\\
$(1,2,1)$&$2\mathbf1_{\{1\}}+P=0$&$1+u+u^2+2y^2+x+2xu+xy$\\
$(2,0,1)$&$\mathbf1_{\{1\}}+P+\lambda_a=0$&$1+u+u^2+y+2yu+xy+2x^2$\\
$(2,2,2)$&$P+\lambda_a=0$&$2+u+2y+2yu+2y^2+2x+x^2$.
\end{tabular}
\end{center}
The verifier checks each assertion at all $81$ group elements.
Here $|R|=2s+1$, so every displayed quadratic has moment zero by
Lemma~\ref{lem:fh-local-moments}. Moreover
$\sum_{g\in J}f_R(g)\lambda_a(g)=C_{u_a}=t$ by exact convolution
of the original projected factors. Summing the four identities
against $f_R$ gives $0=1,2,2,2$, respectively. All four kernels
are therefore excluded. This completes the zero quotient branch.
\end{proof}

\begin{proof}[Proof of Theorem~\ref{thm:fh-small-centre}]
The case $Q=\varnothing$ is Lemma~\ref{lem:fh-qzero}.
For $|Q|=1$, Lemma~\ref{lem:fh-k8} excludes K8 and
Lemma~\ref{lem:fh-qone-K3} excludes K3. Lemma~\ref{lem:fh-exceptions}
and the weighted identity then supply a single functional on $W$
for all extensions, including the central ones.
\end{proof}

\begin{remark}
Each use of an extension in this proof is individual. A union of
three finite product sets is used only to select one of three
separately admissible augmented kernels. No list containing several
extensions simultaneously is assumed to be product-one-free.
Likewise, positions in fibres of weight zero are never discarded;
they remain in the original complement in the moment arguments.
\end{remark}

%% file: sections/q1.tex
\section{A finite obstruction and the cases with at most one central term}
\label{sec:q1}

We give the finite input and its transfer to actual group sequences explicitly.
Throughout this section, subsets are subsets of positions; equal group elements
at different positions remain distinct.

Put $P=\mathbb F_3^4$, with
\[
 \omega(x,y)=x_0y_1-x_1y_0+x_2y_3-x_3y_2.
\]
Let $\mathcal A$ consist of the unordered pairs $\{x,y\}$ of nonzero
vectors with $\omega(x,y)\ne0$. A \emph{balanced triple} is an unordered
triple $T=\{x,y,z\}$ such that
\[
 \omega(x,y)=\omega(y,z)=\omega(z,x)\ne0
\]
in either cyclic ordering. Let $\mathcal D$ be the balanced triples with
$x+y+z=0$, and let $\mathcal I$ be the remaining balanced triples.
There are $2160$ edges, $720$ triples in $\mathcal D$, and $5760$ in
$\mathcal I$. The triples in $\mathcal I$ are linearly independent:
if a balanced triple lies in a plane, writing $z=ax+by$ forces
$a=b=-1$.

For $T\in\mathcal I$ write $s_T=\sum_{x\in T}x$ and define
\[
 L(T)=T\cup\{s_T,-s_T\}.
\]
These are five distinct points. Indeed, every corner is perpendicular to
$s_T$, whereas the three corners are independent. We call $T,U$
\emph{incompatible} when $L(T)\cap L(U)=\varnothing$.
Equivalently, their sum lines differ, they share no exact corner value,
and no corner of either lies on the sum line of the other.
In the last alternative the two sum lines are automatically perpendicular.
This equivalence specifies the compatibility relation without choices of
representatives for projective points.

A \emph{stress} is a family $(\lambda_T)_{T\in\mathcal D\cup\mathcal I}$
in $\mathbb F_3$ satisfying
\begin{equation}\label{eq:q1-stress}
 \sum_{T\supset e}\lambda_T=0\qquad(e\in\mathcal A).
\end{equation}
Its charge is $\ell(\lambda)=\sum_T\lambda_T$.

\begin{theorem}[Finite stress obstruction]\label{thm:finite-q1f}
Every stress of nonzero charge has three pairwise incompatible triples in
its nonzero support in $\mathcal I$.
\end{theorem}

\begin{proof}[Computer-assisted proof]
We describe the exhaustive certificate and the checks that make it a proof.
For a function $\phi:\mathcal A\to\mathbb F_3$, write
$\phi(T)=\sum_{e\subset T}\phi(e)$.
The ancillary certificate supplies a potential $U$ satisfying
\begin{equation}\label{eq:q1-universal}
 U(D)=0\quad(D\in\mathcal D),\qquad U(T)=1\quad(T\in\mathcal I).
\end{equation}
Its $2160$ entries are checked by direct addition on all $6480$ triples.
Call a family $H\subseteq\mathcal I$ \emph{good} if there is a potential
$\psi$ with
\begin{equation}\label{eq:q1-good}
 \psi(D)=1\quad(D\in\mathcal D),\qquad \psi(T)=0\quad(T\in H).
\end{equation}
If the independent support of a stress is contained in a good family,
pairing \eqref{eq:q1-stress} with $\psi+U$ gives
\[
 0=\sum_{e\in\mathcal A}(\psi+U)(e)\sum_{T\supset e}\lambda_T
   =\sum_T\lambda_T=\ell(\lambda).
\]
It therefore suffices to show that every family containing no three
pairwise incompatible members is good.

All geometry is enumerated from the displayed alternating form.
The symplectic similitude group has $103680$ elements, is transitive on
$\mathcal I$, and preserves incompatibility. Fix the first triple $T_0$
in lexicographic order under the point encoding
$x_0+3x_1+9x_2+27x_3$. Its stabilizer has order $18$ and has $244$ orbits
on the $4284$ triples incompatible with $T_0$.
Explicit linear similitudes identify the ordered-pair orbit representatives
under interchange of the two members. This involution has $34$ fixed
points and $105$ pairs; hence $139$ representative roots suffice.
The ancillary verification checks additivity, bijectivity, the similitude
identity, and the claimed pair images for every supplied map. It also
recomputes the group action and all orbit coverage.

For each representative incompatible pair $a,b$, a binary certificate
tree treats every family $F$ with $\{a,b\}\subseteq F$ and no incompatible
triple. A node maintains $S\subseteq F\subseteq H$, initially
\[
 S=\{a,b\},\qquad
 H=\mathcal I\setminus(\operatorname{Inc}(a)\cap
                         \operatorname{Inc}(b)),
\]
where $\operatorname{Inc}(a)$ is the set of triples incompatible with $a$.
At a branch element $v\in H\setminus S$, the exclusion child removes $v$
from $H$. The inclusion child adjoins $v$ to $S$ and removes from $H$
every element in
\[
 \bigcup_{s\in S\cap\operatorname{Inc}(v)}
       \bigl(\operatorname{Inc}(s)\cap\operatorname{Inc}(v)\bigr).
\]
Every allowed family belongs to one of these two children, so the
branching is exhaustive. At every leaf, the verifier recomputes a
solution of \eqref{eq:q1-good} for the complete remaining hull and then
checks its triangle sums directly. Thus correctness does not depend on
the search program or its linear algebra.

The $139$ trees have $15793967$ internal nodes and $15794106$ leaves.
The consumer checks the complete tree grammar, every branch index,
the root hulls, both child hulls, and all $15794106$ leaf potentials.
It rejects a truncated tree. All calculations are over $\mathbb F_3$.
The dependent triples partition $\mathcal A$: the unique dependent
triple containing $\{x,y\}$ is $\{x,y,-x-y\}$. This gives $1440$ free
edge variables after imposing the dependent equations in
\eqref{eq:q1-good}; nevertheless every resulting potential is also
checked in the original $2160$ edge coordinates.

This covers every family containing an incompatible pair. A nonempty
pairwise compatible family can be enlarged by a triple incompatible
with one of its members; the enlarged family still contains no
incompatible triple. The empty family is good by the dependent-edge
partition. Goodness passes to subfamilies, completing the proof.

The ancillary package contains the potentials, the explicit orbit maps,
the complete tree protocols, and their independent consumers. The
reproduction command and the checksums are recorded in its README.
\end{proof}

We now prove the transfer, including the treatment of repeated projected
values and central positions lying in the commutator subgroup.

Write $G_r=E_2\times C_3^r$ in coordinates $(v,z)$, where
$V=\mathbb F_3^{r+4}$ carries an alternating form of rank four and
\[
 (v,z)(w,t)=(v+w,z+t+2\omega(v,w)).
\]
Let $\pi:V\to P$ be the quotient by the radical. A position $(v_i,z_i)$
is central precisely when $\pi(v_i)=0$. Put $d=r+4$.

\begin{lemma}[Product fibres of zero-sum subsets]\label{lem:q1-fibres}
In a product-one-free sequence in $G_r$, every nonempty positional subset
$B$ satisfying $\sum_{i\in B}v_i=0$ has one of the following forms:
\begin{enumerate}
\item all its positions commute, and $\sum_{i\in B}z_i\ne0$;
\item its graph of noncommuting positions is a balanced triangle and isolated
vertices, and $\sum_{i\in B}z_i=0$.
\end{enumerate}
\end{lemma}

\begin{proof}
For an ordering $\rho$ the correction to the sum of the heights is
$2\sum_{a<b}\omega(v_{\rho(a)},v_{\rho(b)})$.
For three positions, writing the cyclic pairings as $a,b,c$, its values
are the signed values of
\[
 2(a+b-c),\quad 2(a-b+c),\quad 2(-a+b+c).
\]
They form $\{0\}$ if all pairings vanish, $\{1,2\}$ if there is one
edge or a balanced triangle, and $\mathbb F_3$ otherwise. A triple of
the last kind gives the full fibre for any containing subset: place
the triple consecutively and vary its internal order. Similarly two
disjoint balanced triangles give the full fibre, since
$\{1,2\}+\{1,2\}=\mathbb F_3$.

In the absence of a full triple, every path of length two must close
to a balanced triangle. Thus each nontrivial connected component of
the nonorthogonality graph is a clique. Four vertices cannot have all
their triples balanced: balance on $123$ and $124$ gives equal
pairings of $1$ with $3$ and $4$, whereas balance on $134$ requires
them to be opposite nonzero values. The components therefore have
size at most three. If the vector sum is zero, every row sum of the
Gram matrix is zero, so an edge component is impossible. There can
be at most one triangle component. The product fibres of the two
remaining types are respectively
$\{\sum z_i\}$ and $\sum z_i+\{1,2\}$; avoiding zero gives the stated
height conditions.
\end{proof}

\begin{lemma}[The actual stress]\label{lem:q1-actual-stress}
A product-one-free sequence of length $2d+3$ in $G_r$ determines a
stress on $\mathcal D\cup\mathcal I$ of charge one. Every triple in
its nonzero support has an actual balanced zero-sum completion in
the original sequence.
\end{lemma}

\begin{proof}
Let $\mathcal B$ be all positional subsets of the second type in
Lemma~\ref{lem:q1-fibres}, and write $\epsilon(B)=(-1)^{|B|}$.
For an $m$-dimensional space over $\mathbb F_3$, its group algebra is
\[
 \mathbb F_3[\mathbb F_3^m]
 \simeq\mathbb F_3[Y_1,\ldots,Y_m]/(Y_1^3,\ldots,Y_m^3).
\]
The $(2m+1)$st power of its augmentation ideal is zero.
Apply this to the additive lifts $(v_i,z_i)$ in dimension $d+1$.
The product of their $2d+3$ augmentation factors is zero. Its identity
coefficient and Lemma~\ref{lem:q1-fibres} give
\begin{equation}\label{eq:q1-charge-count}
 \sum_{B\in\mathcal B}\epsilon(B)=-1.
\end{equation}
For each indexed nonorthogonal pair $i,j$, apply the same argument in
$\mathbb F_3[V]$ to the $2d+1$ factors outside the pair, taking the
coefficient at $-v_i-v_j$. Every zero-sum subset containing the pair
belongs to $\mathcal B$, whence
\begin{equation}\label{eq:q1-edge-count}
 \sum_{B\in\mathcal B:\,i,j\in B}\epsilon(B)=0.
\end{equation}

Each $B\in\mathcal B$ has a unique indexed balanced core. Its three
projected values are distinct and form a balanced value triple $T(B)$.
Set
\[
 \lambda_T=-\sum_{B\in\mathcal B:\,T(B)=T}\epsilon(B).
\]
Equation~\eqref{eq:q1-charge-count} gives $\ell(\lambda)=1$.
For a fixed nonorthogonal value edge $e=\{x,y\}$, sum
\eqref{eq:q1-edge-count} over all indexed pairs projecting to $e$.
A balanced subset contributes exactly once when its core contains
$e$, and zero otherwise. Indeed, only its three core pairs are
nonorthogonal. This gives \eqref{eq:q1-stress} even when projected
values repeat in the original sequence. Finally, a nonzero aggregate
$\lambda_T$ contains at least one actual completion in its defining sum.
\end{proof}

\begin{theorem}[Zero or one central position]\label{thm:q01-exclusion}
For every $r\ge0$, a sequence of length $2r+11$ in $G_r$ with at most
one central position has a nonempty product-one subsequence.
\end{theorem}

\begin{proof}
Suppose otherwise and take the stress in Lemma~\ref{lem:q1-actual-stress}.
For an independent balanced core with projected triple $T$, the
orthogonal complement of the span of $T$ in $P$ is the line
$\langle s_T\rangle$. Thus every noncentral isolate in any completion
projects onto that line.

Choose actual completions $B_T,B_U$ of two incompatible supported
triples. They cannot be disjoint: their union would be vector-zero
and contain two disjoint balanced triangles, hence have full product
fibre. If they share a noncentral position, a core--core intersection
gives a common exact corner, a core--isolate intersection puts a corner
on the other sum line, and an isolate--isolate intersection identifies
the sum lines. Each makes $T,U$ compatible. Consequently incompatible
completions can meet only at central positions.

With no central position, no incompatible pair is possible. With a
single central position $g$, completions of an incompatible pair
intersect exactly in $\{g\}$. If there were three pairwise incompatible
supported triples, choose their completions $B_1,B_2,B_3$. Their
remainders $J_i=B_i\setminus\{g\}$ are pairwise disjoint, each has
vector sum $-v_g$, and their union is vector-zero because $3v_g=0$.
It contains three disjoint balanced cores and has full product fibre,
again a contradiction. This argument includes $v_g=0$.

The nonzero independent support therefore contains no incompatible
triple. Theorem~\ref{thm:finite-q1f} forces charge zero, contrary to
Lemma~\ref{lem:q1-actual-stress}.
\end{proof}

%% file: sections/compression.tex
% Compression and its consequences using the exclusions established above.
\section{Compression of central blocks}
\label{sec:central-compression}

Throughout this section, lists are indexed: equal group elements at
different positions remain different positions.  A list is
\emph{product-one-free} if no nonempty selection of its positions can be
ordered to have product the identity.  For a finite group $G$, write
\[
 N_A=\sum_{a\in A}X^a\in\mathbb F_3[G]
 \qquad(A\leq G),
\]
where $X^gX^{g'}=X^{gg'}$.
We use additive notation inside elementary abelian central subgroups.
For a central list $Q=(q_1,\ldots,q_m)$ put
\[
 F_Q=\prod_{i=1}^m(1-X^{q_i}).
\]
The factors in this product commute.  The arguments below concern actual
products of selections of indexed original positions, in addition to
identities in the group algebra.

\begin{lemma}[Degree of a central coefficient function]
\label{lem:central-coefficient-degree}
Let $Z\cong\mathbb F_3^h$ and let $Q$ be a list of $m$ elements of $Z$.
The coefficient function
\[
 f_Q(z)=[X^z]F_Q\qquad(z\in Z)
\]
is represented by a polynomial of total degree at most $2h-m$, with the
interpretation that it is zero if $m>2h$.  If $Q$ is product-one-free,
then $f_Q(0)=1$, and consequently $m\leq2h$.
\end{lemma}

\begin{proof}
Choose coordinates $z=(z_1,\ldots,z_h)$ on $Z$.  Its identity indicator is
\[
 \delta_0(z)=\prod_{j=1}^h(1-z_j^2).
\]
For $q\in Z$, set $\Delta_qf(z)=f(z)-f(z-q)$.  Expanding the product
defining $F_Q$ gives
\[
 f_Q=\Delta_{q_1}\cdots\Delta_{q_m}\delta_0.
\]
A finite difference lowers total degree by at least one, unless the
result is zero.  This proves the degree assertion.  If $Q$ is
product-one-free, only its empty selection has sum zero, so its signed
coefficient at zero is $1$.  This also proves the final assertion.
\end{proof}

\begin{lemma}[Odd central factorization]
\label{lem:odd-central-factorization}
Let $Q$ be a product-one-free central list spanning an elementary
abelian subgroup $Z\cong\mathbb F_3^h$, where $h\geq1$, and suppose
$|Q|=2h-1$.  Put $w=\sum_{q\in Q}q$.  There is a nonzero linear form
$\ell:Z\longrightarrow\mathbb F_3$ satisfying $\ell(w)=1$, and, for
$U=\ker\ell$,
\begin{equation}
 \label{eq:odd-central-factorization}
 F_Q=N_U(1-X^w),\qquad \dim U=h-1.
\end{equation}
Every element of $U\cup(w+U)$ is the sum of a selection of positions of
$Q$; the selection is allowed to be empty.
\end{lemma}

\begin{proof}
Lemma~\ref{lem:central-coefficient-degree} gives
$f_Q(z)=1+\ell(z)$ for a linear form $\ell$.
The full selection is the only selection with sum $w$: any other such
selection would have a nonempty complementary selection of sum zero.
Thus
\[
 f_Q(w)=(-1)^{2h-1}=-1,
 \qquad \ell(w)=1.
\]
In particular, $\ell\ne0$.  The three cosets of $U$ in $Z$ are
$U,w+U,2w+U$.  On these cosets $1+\ell$ has respective values
$1,-1,0$.  These are precisely the coefficients of
$N_U(1-X^w)$, proving~\eqref{eq:odd-central-factorization}.
A nonzero coefficient is a nonzero signed sum of selections with the
corresponding sum, so at least one such selection exists.
\end{proof}

\begin{theorem}[Compression of an odd central block]
\label{thm:central-compression}
Let $L$ be a product-one-free list in a finite group $G$.  Suppose a
sublist $Q$ of its central positions spans an elementary abelian
$3$-subgroup $Z$ of dimension $h\geq1$ and has length $2h-1$.
Let $w$ and $U$ be as in Lemma~\ref{lem:odd-central-factorization}, and
let $\pi:G\longrightarrow G/U$ be the quotient map.  Delete all the
positions of $Q$ and introduce one new position with value $\pi(w)$:
\[
 L_* = \pi(L\setminus Q)\mathbin{\sqcup}(\pi(w)).
\]
Then $L_*$ is product-one-free and
\begin{equation}
 \label{eq:central-compressed-length}
 |L_*|=|L|-2(h-1).
\end{equation}
Moreover, if
\[
 \mathcal E_G(L)=\{g\in G:L\mathbin{\sqcup}(g)
                         \text{ is product-one-free}\},
\]
then
\begin{equation}
 \label{eq:extension-image-inclusion}
 \pi(\mathcal E_G(L))\subseteq\mathcal E_{G/U}(L_*).
\end{equation}
The new position is a block position, not a relabelled single original
position.
\end{theorem}

\begin{proof}
Suppose first that a nonempty product-one word in $L_*$ does not use the
new position.  The corresponding word in positions of $L\setminus Q$
has a product $p\in U$.  By Lemma~\ref{lem:odd-central-factorization},
a selection from the disjoint block $Q$ has product $p^{-1}\in U$.
Concatenating these words gives a nonempty product-one word in $L$.

If the new position is used, its centrality allows it to be moved to
the end of the word.  The remaining original positions have a product
$p$ satisfying $pw\in U$, and hence $p^{-1}\in w+U$.  Again a selection
from $Q$ has product $p^{-1}$.  The remaining word cannot be empty,
since $w\notin U$ and therefore $\pi(w)\ne1$.  Concatenation again
contradicts product-one-freeness of $L$.

The length identity counts the deleted $2h-1$ positions and the one new
position.  Finally, apply the same argument to
$L\mathbin{\sqcup}(g)$ for each $g\in\mathcal E_G(L)$, retaining the
added position outside $Q$.  Its compression is exactly
$L_*\mathbin{\sqcup}(\pi(g))$.  This proves the inclusion.
\end{proof}

\begin{remark}[The implication is one-way]
\label{rem:compression-not-equivalence}
Equation~\eqref{eq:extension-image-inclusion} is an inclusion, not an
equality, and product-one-freeness of the compressed list does not
imply product-one-freeness of the original list.  For example, let
$e,f$ be a basis of $\mathbb F_3^2$ and set
\[
 Q=(e,f,e+f),\qquad w=2e+2f,\qquad U=\langle e-f\rangle.
\]
The list $Q$ is product-one-free and satisfies
$F_Q=N_U(1-X^w)$.  Its extension by $w$ is not product-one-free, since
$(e+f)+w=0$.  In contrast, the compressed extension
$(\pi(w),\pi(w))$ is product-one-free in $\mathbb F_3^2/U$.
Vanishing group-algebra coefficients do not assert that the
corresponding products are absent from the original list.
\end{remark}

\begin{proposition}[Saturated central blocks]
\label{prop:central-saturation}
Let $L$ be product-one-free and let a sublist $Q$ of its central
positions span an elementary abelian $3$-subgroup $Z$ of dimension
$h$, with $|Q|=2h$.  Then
\begin{equation}
 \label{eq:central-saturation}
 F_Q=N_Z.
\end{equation}
For $\pi:G\longrightarrow G/Z$, the list
$L_* = \pi(L\setminus Q)$ is product-one-free, has length $|L|-2h$,
and satisfies
$\pi(\mathcal E_G(L))\subseteq\mathcal E_{G/Z}(L_*)$.
\end{proposition}

\begin{proof}
Lemma~\ref{lem:central-coefficient-degree} makes the coefficient
function constant, with value $1$ at the identity, proving the
factorization.  In particular, every element of $Z$ is a product of
a selection from $Q$.  A nonempty identity word in the quotient would
lift to an original word with product in $Z$; a disjoint selection
from $Q$ cancels that product.  The same proof applies after any
single admissible extension.  The assertion also covers $h=0$ and an
empty block.
\end{proof}

\begin{proposition}[Transfer of an ordered norm factorization]
\label{prop:central-norm-transfer}
In Theorem~\ref{thm:central-compression}, any ordering of the augmentation
factors $1-X^g$ indexed by $L_*$ whose product is $N_{G/U}$ lifts to an
ordering of all the original augmentation factors indexed by $L$ whose
product is $N_G$.  The corresponding
statement holds for Proposition~\ref{prop:central-saturation}.
\end{proposition}

\begin{proof}
For a central subgroup $U$, define the linear map
\[
 I_U:\mathbb F_3[G/U]\longrightarrow\mathbb F_3[G],
 \qquad I_U(X^{\pi(g)})=N_UX^g.
\]
This is well defined and injective, since different cosets have
disjoint supports.  Expanding an ordered product gives
\begin{equation}
 \label{eq:norm-induction-product}
 I_U\left(\prod_j(1-X^{\pi(g_j)})\right)
       =N_U\prod_j(1-X^{g_j}),
 \qquad I_U(N_{G/U})=N_G.
\end{equation}
In the odd case, lift the new factor using $w$.
Equation~\eqref{eq:odd-central-factorization} then replaces
$N_U(1-X^w)$ by the product of all the original factors from $Q$.
Since these factors are central, they may be inserted consecutively at
the position of the new factor.  Every original position is used
exactly once.  In the saturated case, use
\eqref{eq:central-saturation} instead.
Only linearity and the explicitly expanded
identity~\eqref{eq:norm-induction-product} are used; $I_U$ is not
asserted to be an algebra homomorphism.
\end{proof}

\subsection{Central quotients and critical lists}

For clarity, the needed quotient structure is recorded explicitly.
Let $K_m$ denote the extraspecial group of exponent $3$ and order
$3^{2m+1}$, so $K_1=H_{27}$ and $K_2=E_2$.  Write
$G=K_m\times C_3^a$ and $C=G'$, a subgroup of order $3$.
If $U\leq Z(G)$ has dimension $k$, then
\begin{equation}
 \label{eq:central-quotient-type}
 G/U\cong
 \begin{cases}
 K_m\times C_3^{a-k}, & U\cap C=1,\\
 C_3^{2m+1+a-k}, & C\leq U.
 \end{cases}
\end{equation}
In the first case $k\leq a$.  Indeed, the extraspecial factor embeds
in the quotient, since its intersection with the central subgroup
$U$ lies in $C$.  A central complement to its commutator subgroup in
$Z(G)/U$ supplies the displayed direct factor.  In the second case
the quotient is abelian of exponent $3$, with the displayed order.
These cases exhaust the possibilities because $C$ has order $3$.
Also, in the first case no original noncentral element becomes
central: if $[g,x]\ne1$, then $[g,x]\in C\setminus U$ remains
nontrivial in the quotient.

\begin{corollary}[A central deficit condition]
\label{cor:central-deficit}
Let $S$ be a product-one-free list of length $2r+11$ in
$G=E_2\times C_3^r$.  Let $Q$ be its entire sublist of positions in
$Z(G)$, and let $h=\dim_{\mathbb F_3}\langle Q\rangle$ in the full
centre.  If $h\geq1$, then
\begin{equation}
 \label{eq:critical-central-deficit}
 |Q|\leq2h-2.
\end{equation}
\end{corollary}

\begin{proof}
The abelian bound in Lemma~\ref{lem:central-coefficient-degree} gives
$|Q|\leq2h$.  Suppose first that $|Q|=2h-1$, and apply
Theorem~\ref{thm:central-compression}, with $\dim U=h-1$.

If $U\cap G'=1$, equation~\eqref{eq:central-quotient-type} gives
$G/U\cong E_2\times C_3^{r-h+1}$.  The compressed length is
$2(r-h+1)+11$.  No noncentral position becomes central, and all the
original central positions have been replaced by the one new position.
Thus the compressed list has exactly one central position, contrary
to Theorem~\ref{thm:q01-exclusion}.

If $G'\leq U$, the quotient is elementary abelian of dimension
$d=r+6-h$.  The compressed list has length
$2r+13-2h=2d+1$, contrary to
Lemma~\ref{lem:central-coefficient-degree}.

It remains to exclude $|Q|=2h$.  Apply
Proposition~\ref{prop:central-saturation}.  If $\langle Q\rangle$
meets $G'$ trivially, the resulting list is a central-free critical
list in $E_2\times C_3^{r-h}$. Such lists are excluded by
Theorem~\ref{thm:q01-exclusion}.
If $G'\leq\langle Q\rangle$, the quotient has dimension
$d=r+5-h$ and the resulting list again has length $2d+1$, which is
impossible.  This proves~\eqref{eq:critical-central-deficit}.
\end{proof}

\begin{remark}
When $h=1$, the odd compression has $U=1$ and does not shorten the list;
this endpoint uses the one-central-position exclusion directly.
For $h>1$, the same reduction treats arbitrarily large odd
central blocks.  The rank $h$ here is the rank of the original central
values in the full centre, not their rank after quotienting by $G'$.
\end{remark}

\subsection{Lifting extension hyperplanes}

The following consequence uses the theorem for at most one central
position established in Theorem~\ref{thm:fh-small-centre}.
All linear forms refer to the full Baker--Campbell--Hausdorff vector
space of the exponent-$3$, class-$2$ group, including its commutator
coordinate.

\begin{corollary}[Extension hyperplanes at central deficit zero or one]
\label{cor:fh-central-deficit}
Let $L$ be a product-one-free list of length $2s+5$ in
$H=H_{27}\times C_3^s$.  Let $Q$ be its entire sublist of
$H$-central positions, and put $h=\dim_{\mathbb F_3}\langle Q\rangle$.
If $|Q|\geq2h-1$, there is a homogeneous linear form $\nu$ on the
full BCH vector space of $H$ such that
\[
 \nu(g)=1\qquad\text{for every }g\in\mathcal E_H(L).
\]
This includes central admissible extensions.
\end{corollary}

\begin{proof}
When $h=0$, the list $Q$ is empty and the assertion is
Theorem~\ref{thm:fh-small-centre}.  For $h\geq1$, the abelian bound
leaves exactly the cases $|Q|=2h-1$ and $|Q|=2h$.

In the odd case, compress by $U$ of dimension $h-1$ as in
Theorem~\ref{thm:central-compression}.  If $H'$ survives, the quotient
is $H_{27}\times C_3^{s-h+1}$ and the compressed list has length
$2(s-h+1)+5$ with exactly one central position.  Apply
Theorem~\ref{thm:fh-small-centre} in that quotient and pull its linear
form back along $\pi$.  The quotient map is linear on the full BCH
spaces, and~\eqref{eq:extension-image-inclusion} shows that the pulled
back form has value $1$ on every original admissible extension.

If $H'\leq U$, the quotient is elementary abelian of dimension
$d=s+4-h$ and the compressed length is $2d-1$.  Its group-algebra
coefficient function is $1+\lambda$ by
Lemma~\ref{lem:central-coefficient-degree}.  For an admissible
extension $a$, no selection of the compressed list can sum to $-a$;
otherwise that selection, together with the new position $a$, would
sum to zero.  Therefore
\[
 0=f_{L_*}(-a)=1-\lambda(a),
\]
and $\lambda(a)=1$.  Pulling back $\lambda$ proves the assertion.

In the saturated case, use Proposition~\ref{prop:central-saturation}
instead.  If $H'$ survives, the quotient is
$H_{27}\times C_3^{s-h}$ and the compressed list has length
$2(s-h)+5$ with no central positions; apply
Theorem~\ref{thm:fh-small-centre} again.  If $H'$ is killed, the
quotient is elementary abelian of dimension $d=s+3-h$ and the
compressed length is again $2d-1$, so the preceding affine
coefficient argument applies unchanged.
\end{proof}